\documentclass[a4paper,12pt]{article}

\usepackage{amssymb}
\usepackage{graphicx}
\usepackage{amsmath}
\usepackage{amsfonts}

\newcommand{\singlespacing}{\let\CS=\@currsize\renewcommand{\baselinestreatch}{1.0}\tiny\CS}
\newcommand{\doublespacing}{\let\CS=\@currsize\renewcommand{\baselinestreatch}{1.5}\tiny\CS }

\newtheorem{theorem}{Theorem}[section]
\newtheorem{corollary}{Corollary}[section]

\newtheorem{lemma}{Lemma}[section]

\numberwithin{equation}{section}

\begin{document}
\begin{center}
\textbf{\Large {The Critical Point Equation And Contact Geometry}}
\end{center}
\centerline{Amalendu Ghosh$^1$ and Dhriti Sundar Patra$^2$}

\newtheorem{Theorem}{\quad Theorem}[section]
\newtheorem{Definition}[Theorem]{\quad Definition}
\newtheorem{Corollary}[Theorem]{\quad Corollary}
\newtheorem{Lemma}[Theorem]{\quad Lemma}
\newtheorem{Example}[Theorem]{\emph{Example}}
\newtheorem{Proposition}[Theorem]{Proposition}
\numberwithin{equation}{section}
\noindent\\
\textbf{Abstract:} {In this paper, we consider the CPE conjecture in the frame-work of $K$-contact and $(\kappa, \mu)$-contact manifolds. First, we prove that if a complete $K$-contact metric satisfies the CPE is Einstein and is isometric to a unit sphere $S^{2n+1}$}. Next, we prove that if a non-Sasakian $ (\kappa, \mu) $-contact metric satisfies the CPE, then $ M^{3} $ is flat and for $ n > 1 $, $ M^{2n+1} $ is locally isometric to $ E^{n+1}\times S^{n} (4) $.

\noindent\\
\textbf{Mathematics Subject Classification 2010}: 53C25, 53C20, 53C15

\noindent\\
\textbf{Keywords}: Total scalar curvature functional, Critical point equation, K-contact manifold, $(\kappa, \mu)$-contact manifold, Einstein manifold.

\section{Introduction}
Let $\mathcal{M}$ denote the set of Riemannian metrics on a compact orientable manifold $M^n$ of unit volume. Given a Riemannian metric $g\in \mathcal{M}$, the total scalar curvature functional $\mathcal{S} :\mathcal{M} \longrightarrow  R$ is defined by
\begin{eqnarray*}
\mathcal{S}(g) =\int_{M}r_{g}dv_{g},
\end{eqnarray*}
where $r_{g}$ is the scalar curvature and $dv_{g}$ the volume form determined by the metric and orientation. The functional $\mathcal{S}$ restricted over $\mathcal{M}$ is known as  Einstein-Hilbert functional and its critical points are the Einstein metrics ( see Chapter $2$ in \cite{AB}).
In \cite{C}, Corvino proved that $\lambda$ is a nontrivial solution of $\mathcal{L}^{*}_g(\lambda) = 0$ if and only if the warped product metric $g^* = g - \lambda^2 dt^2$ is Einstein. Here, $\mathcal{L}^{*}_g(\lambda)$ is the formal $L^2$-adjoint of the linearized scalar curvature operator $\mathcal{L}_g(\lambda)$ and is defined as
\begin{eqnarray}\label{1.1}
 \mathcal{L}^{*}_g(\lambda) = -(\Delta_g\lambda)g + Hess_{g} \lambda - \lambda Ric_{g},
\end{eqnarray}
where $\Delta_g$, $Ric$ and $Hess\lambda$ are respectively the Laplacian, the Ricci tensor and the Hessian of the smooth function $\lambda$ on $M$.\\

The classical Yamabe problem says that any compact manifold carries many smooth Riemannian metrics with constant scalar curvature. So we may introduce the set of constant scalar curvature-metrics as follows:
$$\mathcal{C} =\{g \in M|r_{g} = constant\}.$$
The Euler-Lagrange equation of Hilbert-Einstein functional restricted to $\mathcal{C}$ on a given compact oriented manifold $(M, g)$ can be written as the following critical point equation (shortly, CPE)
\begin{eqnarray}\label{1.2}
  \mathcal{L}^{*}_g(\lambda)= {Ric}^{o}_{g},
\end{eqnarray}
where ${Ric}^{o}_{g}$ denotes the traceless Ricci tensor of $M$. The function $\lambda $ is known as the potential function. It is interesting to point out that if $\lambda$ is constant in the equation (\ref{1.1}), then $\lambda = 0$ and $g$ becomes Einstein. Therefore, from now on, we consider a metric $g$ with a non-trivial potential function $\lambda$ as a solution of CPE and is denoted by $(g,\lambda)$.
Using (\ref{1.1}), one can express the equation (\ref{1.2}) in the following form
\begin{eqnarray}\label{1.3}
Hess_{g} \lambda + (\frac{r}{n-1}g - Ric_{g} )\lambda = Ric_{g} - \frac{r}{n}g .
\end{eqnarray}
Thus, we consider the following definition
\begin{Definition}
A compact oriented Riemannian manifold $(M, g)$ of dimension $n\geq 3$ with constant scalar curvature and volume $1$ together with a non-constant smooth potential function $\lambda$ satisfying the equation (\ref{1.2}) is called a Critical Point Equation metric.
\end{Definition}
In \cite{AB}, A. Besse conjectured that the solution of the CPE is Einstein (see \cite{AB}, p. 128). Since then, several mathematician have attempted to prove the CPE conjecture. However, the conjecture is yet to be proved. Although, some partial answers were obtained under some curvature assumptions. For example, Lafontaine proved that the CPE conjecture is true under conformally flat assumption with $Ker\mathcal{L}^{*}_g(\lambda) \ne 0$. Recently, Barros and Ribeiro Jr \cite{ABC} proved that the CPE conjecture is also true for half conformally flat.
Another partial proof of the CPE conjecture was presented by Yun, Chang and Hwang \cite{YCH}. They proved that if $(g,\lambda)$ is a non-trivial solution of the CPE on an $n$-dimensional compact Riemannian manifold $M$ and satisfies one of the following conditions $(i)$ Ricci tensor of $g$ is parallel $(ii)$ $g$ has harmonic curvature or $(iii)$ $g$ is conformally flat, then $(M,g)$ is isometric to a standard sphere. In \cite{HC}, Hwang proved that the CPE conjecture is also true under certain conditions on the bounds of the potential function $\lambda$. Very recently, Nato \cite{Na} deduced a necessary and sufficient condition on the norm of the gradient of the potential function for a CPE metric to be Einstein. \\

In this paper, we consider the CPE conjecture in the frame-work of $K$-contact manifolds and $ (\kappa, \mu) $-contact manifolds. Let $M$ be a $(2n+1)$-dimensional contact manifold with $(\varphi, \xi, \eta)$ as its almost contact structure. A Riemannian metric is said to be an associated metric if it satisfies $g(\varphi X, \varphi Y) = g(X, Y) - \eta(X)\eta(Y).$ Then  $M^{2n+1}$ is said to be an almost contact metric manifold with $(\varphi, \xi, \eta, g)$ as its almost contact metric structure. Moreover, if $\xi$ is Killing, then $M$ is said to be $K$-contact (see section 2). This raises the question whether the CPE conjecture is true in the frame-work of $K$-contact manifold and $ (\kappa, \mu) $-contact manifold. In section $3$, we characterize that if a complete $K$-contact metric $g$ satisfies the CPE (\ref{1.2}), then it is Einstein and isometric to a unit sphere $S^{2n+1}$. We also classify $(\kappa, \mu)$-contact metric satisfying the CPE.

\section{Preliminaries}
In this section, we recall some basic definitions and formulas on a contact metric manifold which will be useful for the establishment of our results. A Riemannian manifold of dimension $ (2n + 1)$ is said to be a contact manifold if it admits a global $ 1 $-form $ \eta $ such that $ \eta \wedge ( d\eta )^{n}$ is non-vanishing everywhere on $ M $. This $ 1 $-form is known as contact form. Corresponding to this $ \eta $ one can find a unit vector field $ \xi $, called the Reeb vector field, such that $ \eta ( \xi ) = 1 $ and $ d\eta ( \xi, .) = 0 .$ It is well-known that every contact manifold admits an underlying almost
contact structure $(\varphi, \xi, \eta)$, where $\varphi$ is a global tensor field of type $(1, 1)$, such that
$\hskip 0.1cm \eta(X) = g(X, \xi), \hskip 0.1cm\varphi \xi = 0, \hskip 0.1cm\eta \circ \varphi = 0, \hskip 0.1cm\hskip 0.1cm\varphi^2 = -I + \eta\otimes \xi.$
Further, an almost contact structure is said to be contact metric if it satisfies
$$\hskip 0.1cm d\eta(X,Y)=g(X,\varphi Y), \hskip 0.1cm g(\varphi X, \varphi Y) = g(X,Y) - \eta(X)\eta(Y).$$
A Riemannian manifold $ M^{2n+1} $ together with the structures $( \varphi, \xi, \eta, g )$ is said to be a contact metric manifold.
We now define two operators $ h $ and $ l $ by $ h = \frac{1}{2}\pounds_{\xi}\varphi $ and $ l = R(., \xi )\xi$. These tensors are self-adjoint and satisfy $Tr h = 0 $, $Tr h\varphi =0 $, $l\xi = 0$ and $ h\varphi = - \varphi h $. On a contact metric manifold the following formulas are valid \cite{Blair}
\begin{eqnarray}\label{2.1}
\nabla_{X}\xi = - \varphi X - \varphi h X.
\end{eqnarray}
\begin{eqnarray}\label{2.2}
Ric (\xi, \xi) = g(Q \xi, \xi ) = Tr l = 2n - Tr h^{2}.
\end{eqnarray}
\begin{eqnarray}\label{2.3}
\nabla_{\xi}h = \varphi - \varphi h^2 - \varphi l.
\end{eqnarray}
If the vector field $ \xi $ is Killing (equivalently, $ h = 0 $ or $ Tr l = 2n $), then the contact metric manifold $ M $ is said to be a $ K $-contact. On a $ K $-contact manifold the following formulas are known \cite{Blair}
\begin{eqnarray}\label{2.4}
\nabla_{X}\xi = - \varphi X,
\end{eqnarray}
\begin{eqnarray}\label{2.5}
Q \xi = 2n \xi,
\end{eqnarray}
\begin{eqnarray}\label{2.6}
R(\xi,X)Y=(\nabla_{X}\varphi)Y,
\end{eqnarray}
where $ \nabla $ is the operator of covariant differentiation of $ g $, $Q$ is the Ricci operator associated with the $ (0, 2)$ Ricci tensor $ Ric $ and $ R $ is the Riemann curvature tensor of $ g $.
The following formula also holds for a $K$-contact manifold (as $h = 0$) (see \cite{Blair} p. 116)
\begin{eqnarray}\label{2.7}
(\nabla_{Y}\varphi)X+(\nabla_{\varphi Y}\varphi)\varphi X= 2g(Y,X)\xi-\eta(X)(Y + \eta(Y)\xi),
\end{eqnarray}
 We now deduce some equations which would be used later. Taking covariant differentiation of (\ref{2.7}) along an arbitrary vector field $X$ and using (\ref{2.6}), we get
\begin{equation}\label{2.8}
(\nabla_{X}Q)\xi = Q\varphi X-2n\varphi X.
\end{equation}
As $\xi$ is Killing for a $K$-contact manifold, we have that $\pounds_{\xi}Q = 0.$ Making use of (\ref{2.8}) and (\ref{2.5}), one can easily deduce that
\begin{eqnarray}\label{2.9}
\nabla_{\xi}Q = Q\varphi - \varphi Q.
\end{eqnarray}
A contact metric structure on $M$ is said to be normal if
the almost complex structure on $M\times R$ defined by
$J(X,fd/dt)=(\varphi X-f\xi,\eta(X)d/dt),$
where $f$ is a real function on $M\times R$, is integrable. Equivalently, a contact metric manifold is said to be Sasakian if $(\nabla_{X}\varphi )Y = g(X, Y)\xi - \eta(Y)X,$ or if the curvature tensor satisfies $R ( X, Y)\xi = \eta ( Y )X - \eta ( X )Y.$\\
By a $ (\kappa, \mu) $-contact manifold we mean a contact metric manifold $ M^{2n+1} (\varphi, \xi, \eta, g)$ whose curvature tensor satisfies
\begin{eqnarray}\label{2.10}
R(X, Y)\xi = \kappa\{\eta(Y)X - \eta(X)Y\} + \mu \{\eta(Y)hX - \eta(X)hY \},
\end{eqnarray}
for some real numbers $ (\kappa, \mu).$ This class of contact manifold was introduced by Blair et al. (see \cite{Blair}). In particular, it arises by applying the $ D $-homothetic deformation (\cite{T}):
 $$ \bar{\eta} = a\eta, \bar{\xi}= \frac{1}{a}\xi, \bar{\varphi} = \varphi ,\bar{g} = ag + a(a-1)\eta\otimes\eta,$$
for a positive real constant $a$, to a contact metric manifold satisfying $ R(X, Y)\xi = 0$. Note that $ D $-homothetic deformation preserves Sasakian, K-contact and $(\kappa,\mu)$-contact structures.
It is interesting to point out that the class of $(\kappa,\mu)$-contact structure contains Sasakian manifolds (for $\kappa = 1$) and the trivial sphere bundle $ E^{n+1}\times S^{n}(4) $ (for $\kappa = \mu = 0 $). Examples of non-Sasakian $ (\kappa, \mu) $-contact manifolds are the tangent sphere bundles of Riemannian manifolds of constant curvature $\neq 1$. Further, the equation (\ref{2.7}) determines the curvature completely for $\kappa < 1$. For $ (\kappa, \mu) $-contact manifolds, the following formulas are known (see \cite{Blair})
\begin{eqnarray}\label{2.11}
QX = [2(n - 1) - n\mu ]X + [2(n - 1) + \mu ]hX \nonumber\\
+ [2(1 - n) + n(2\kappa + \mu)]\eta(X)\xi.
\end{eqnarray}
\begin{eqnarray}\label{2.12}
h^{2} = (\kappa - 1)\varphi^{2},
\end{eqnarray}
where $ \kappa \leq 1$.
Moreover, the constant scalar curvature $ r $ of such class is given by
\begin{eqnarray}\label{2.13}
 r = 2n(2(n - 1) + \kappa - n\mu).
\end{eqnarray}

\section{Main Results}
In this section, we consider $K$-contact and $ (\kappa, \mu) $-contact metric satisfying the critical point equation . First, we prove the following:
\begin{lemma}
Let $(g,\lambda)$ be a non-trivial solution of the CPE (\ref{1.2}) on an $n$-dimensional Riemannian manifold $M$.
Then the curvature tensor $R$ can be expressed as
\begin{eqnarray}\label{3.1}
R(X,Y)D\lambda &=& (X\lambda)QY - (Y\lambda)QX + (\lambda +1)(\nabla_{X}Q)Y \nonumber\\
 &-& (\lambda +1)(\nabla_{Y}Q)X + (Xf)Y - (Yf)X.
\end{eqnarray}
\end{lemma}
\textbf{Proof: }Tracing the equation  (\ref{1.3}) implies
$\triangle_{g}\lambda = -\frac{r\lambda}{n-1}$. Thus, the equation (\ref{1.2}) can be exhibited as
\begin{eqnarray}\label{3.2}
\nabla_{X}D\lambda = (\lambda + 1)QX + fX,
\end{eqnarray}
where $f=-r(\frac{\lambda}{n-1} + \frac{1}{n}).$
Taking covariant differentiation of (\ref{3.2}) along an arbitrary vector field Y, we obtain
\begin{eqnarray*}
\nabla_{Y}(\nabla_{X}D\lambda) &=& (Y\lambda)QX  + (\lambda +1)(\nabla_{Y}Q)X + (\lambda +1) Q(\nabla_{Y}X) \nonumber\\
&+& (Yf)X + f\nabla_{Y}X.
\end{eqnarray*}
Repeated application of this equation in the well known expression of the curvature tensor $R(X,Y) = [\nabla_{X},\nabla_{Y}] - \nabla_{[X,Y]},$
gives the required result.

\begin{theorem}
Let $M(\varphi,\xi,\eta,g)$ be a complete K-contact manifold of dimension $(2n+1)$. If $(g,\lambda)$ is a non-constant solution of the critical point equation (\ref{1.2}), then $(M,g)$ is Einstein and isometric to a unit sphere $S^{2n+1}$.
\end{theorem}
\textbf{Proof: }
Replacing $\xi$ instead of $X$ in (\ref{3.1}) and using the equations (\ref{2.5}), (\ref{2.8}) and (\ref{2.9}), we get
\begin{eqnarray}\label{3.3}
R(\xi,Y)D\lambda &=&(\xi\lambda)QY-2n(Y\lambda)\xi- (\lambda +1)\varphi QY \nonumber\\
&+& 2n(\lambda +1) \varphi Y+(\xi f)Y-(Yf)\xi.
\end{eqnarray}
Considering the scalar product of the foregoing equation with an arbitrary vector field $X$ and making use of the equation (2.6) we have
\begin{eqnarray}\label{3.4}
g((\nabla_{Y}\varphi) X,D\lambda)&+&(\xi\lambda)g(QY,X)+2n(\lambda +1) g(\varphi Y,X)\nonumber\\
&-& \{2n(Y\lambda)+(Yf)\}\eta(X)-(\lambda +1)  g(\varphi QY,X)\nonumber\\
&+&(\xi f)g(X,Y)= 0.
\end{eqnarray}
Setting $X=\varphi X, Y=\varphi Y$ in (\ref{3.4}) and adding the resulting equation with (\ref{3.4}) and then using the equation (\ref{2.7}) gives
\begin{eqnarray*}
 2\xi(\lambda+f)g(X,Y)- Y\{(2n+1)\lambda+f\}\eta(X)\\
 -\xi(\lambda+f)\eta(X)\eta(Y)+(\xi\lambda)g(QY,X)+ 4n(\lambda +1)  g(\varphi Y,X)\\
-(\lambda +1)  g(Q\varphi Y+\varphi QY,X) + (\xi\lambda)g(Q\varphi Y,\varphi X) = 0.
\end{eqnarray*}
Anti-symmetrizing the foregoing equation yields
\begin{eqnarray}\label{3.5}
 X\{(2n+1)\lambda+f\}\eta(Y) - Y\{(2n+1)\lambda+f\}\eta(X)\nonumber\\
- 8n(\lambda +1)  g(\varphi X,Y) - 2(\lambda +1)  g(Q\varphi Y+\varphi QY,X)= 0.
\end{eqnarray}
Replacing $X$ by $\varphi X$ and $Y$ by $\varphi Y$ in the preceding equation, we deduce
\begin{equation*}
(\lambda +1)  [g(Q\varphi Y,X)+g(\varphi QY,X)] = 4n(\lambda +1)  g(\varphi Y,X)).
\end{equation*}
Since $\lambda$ is a non-constant smooth function on $M$, the last equation implies that
\begin{eqnarray}\label{3.6}
(Q\varphi + \varphi Q)X = 4n \varphi X,
\end{eqnarray}
for all vector field $X$ in $M$.\\
Let $\{e_{i},\varphi e_{i},\xi \},i= 1,2,3,.....,n$, be a $\varphi-$basis of M such that $Qe_{i} = \rho_{i}e_{i}$. From which, we deduce
$\varphi Qe_{i} = \rho_{i}\varphi e_{i}.$
Substituting $e_{i}$ for Y in the equation (\ref{3.6}) and using the foregoing equation, we obtain $Q\varphi e_{i}  = (4n-\rho_{i}\varphi) e_{i}$. Using the $\varphi $-basis and the equation (\ref{2.5}), the scalar curvature $r$ is given by
$$r = g(Q\xi,\xi) + \sum_{i=1}^{n}[g(Qe_{i},e_{i}) + g(Q\varphi e_{i},\varphi e_{i})] = 2n(2n+1).$$
For a $(2n + 1)$-dimensional $K$-contact manifold, we have $f=-r(\frac{\lambda}{2n} + \frac{1}{2n + 1})$ (follows from lemma (3.1)). Since $r = 2n(2n + 1)$, the last equation reduces to
\begin{equation}\label{3.7}
 (2n+1)\lambda+f = - 2n(constant).
\end{equation}
Now, taking inner product of (\ref{3.3}) with $D\lambda$ and recalling (\ref{3.7}), we get
\begin{equation}\label{3.8}
(\xi\lambda)\{QD\lambda-2nD\lambda\}+(\lambda +1)\{Q\varphi D\lambda-2n\varphi D\lambda\}=0.
\end{equation}
Next, taking $D\lambda$ instead of $Y$ in (\ref{3.6}), we obtain
$$Q\phi D\lambda+\varphi QD\lambda-4n\varphi D\lambda=0.$$
Using the foregoing equation in (\ref{3.8}), we find
\begin{equation}\label{3.9}
(\xi\lambda)\{QD\lambda-2nD\lambda\}+(\lambda +1)\{2n\varphi D\lambda-\varphi QD\lambda\}=0.
\end{equation}
Operating (\ref{3.9}) by $\varphi$ and using (\ref{2.5}) provides
\begin{equation}\label{3.10}
(\xi\lambda)\{\varphi QD\lambda-2n\varphi D\lambda\}+(\lambda +1)\{QD\lambda-2nD\lambda\}=0.
\end{equation}
Equations (\ref{3.10}) and (\ref{3.9}) together imply
\begin{equation}\label{3.11}
\{(\lambda +1) ^2+(\xi\lambda)^2\}(QD\lambda-2nD\lambda)=0.
\end{equation}
If possible, let $(\lambda + 1) ^2+(\xi\lambda)^2=0$ in some open set $\mathcal{O}$ in $M$. Then $\lambda + 1 = 0$ and $\xi\lambda =0$ on $\mathcal{O}$. Since $\lambda$ is not a constant, so $\lambda+ 1 = 0$ is not possible. Consequently, it follows that $QD\lambda-2nD\lambda =0$.
 Now, the covariant differentiation of the foregoing equation along an arbitrary vector field $X$ and then using (\ref{3.2}) gives
\begin{equation*}
(\nabla_{X}Q)D\lambda+\lambda Q^2X+(f-2n\lambda)QX-2nfX=0.
\end{equation*}
Contracting this  over $X$ with respect to an orthonormal field and noting that $r=2n(2n+1)$, we obtain
$|Q|^2=2nr.$
Making use of this and recalling $r=2n(2n+1)$, we compute
\begin{eqnarray*}
&|Q-\frac{r}{2n+1}I|^2 = |Q|^2-\frac{2r^2}{2n+1}+\frac{r^2}{2n+1}=2nr-\frac{r^2}{2n+1}=0.
\end{eqnarray*}
Since the length of the symmetric tensor $Q-\frac{r}{2n+1}I$ vanishes, we must have $Q=\frac{r}{2n+1}I=2nI$.
This shows that $M$ is Einstein with Einstein constant $2n$. Since $M$ is complete, it is compact by Myers' theorem \cite{MSB}. Use of (\ref{3.7}) in (\ref{3.2}) provides
\begin{equation}\label{3.12}
\nabla^2_{g}\lambda = -(\lambda + 2n)g.
\end{equation}
We are now in position to apply Tashiro's Theorem (\cite{YT}): \textit{``If a complete Riemannian manifold $M^n$ of $dim \geq 2$ admits a special concircular field $\rho$ satisfying $\nabla \nabla \rho = (- c^2\rho + b)g$, then it is isometric to a sphere $S^n(c^2)$"} to conclude that $M$ is isometric to a unit sphere $S^{2n+1}$. This completes the proof.

\begin{corollary}
Let $M(\varphi,\xi,\eta,g)$ be a complete and simple connected Sasakian manifold of dimension $(2n+1)$. If $(g,\lambda)$ is a non-constant solution of the critical point equation (\ref{1.2}), then $(M,g)$ is Einstein and isometric to a unit sphere $S^{2n+1}$.
\end{corollary}
\textbf{Proof: } On a Sasakian manifold the Ricci operator $Q$ and $\varphi$ commutes, i.e., $Q\varphi = \varphi Q$ (see \cite{Blair}). Using this in (\ref{3.6}) implies $Q\varphi X = 2n \varphi X$. Substituting X by $\varphi X$ in the last
equation and using (\ref{2.5}) gives $QX = 2nX$. This shows that M is Einstein with Einstein constant $2n$. Rest of the
proof follows from the last Theorem.

\begin{theorem}
Let $ M^{2n+1} (\varphi, \xi, \eta, g)$ be a non-Sasakian $ (\kappa, \mu) $-contact manifold. If $(g,\lambda)$ is a non-constant solution of the critical point equation (\ref{1.2}), then $ M^{3} $ is flat and for $ n > 1 $, $ M^{2n+1} $ is locally isometric to $ E^{n+1}\times S^{n} (4) $.
\end{theorem}
\textbf{Proof:} Taking $ Y = \xi$ in (\ref{2.10}) it follows that
\begin{eqnarray}\label{3.13}
l = -\kappa \varphi^2 + \mu h.
\end{eqnarray}
Using (\ref{2.12}) and (\ref{3.13}) in (\ref{2.3}), we obtain
\begin{eqnarray}\label{3.14}
\nabla_{\xi}h = \mu h\varphi.
\end{eqnarray}
Differentiating covariantly (\ref{2.11}) along $\xi$ and using the equation (\ref{3.14}) yields
\begin{eqnarray}\label{3.15}
(\nabla_{\xi}Q)X = \mu(2(n-1)+ \mu)h\varphi X.
\end{eqnarray}
On the other hand, from (\ref{2.10}) we have $Q\xi = 2n\kappa \xi$.
Differentiating this along an arbitrary vector field X and using (\ref{2.1}) it follows that
\begin{eqnarray}\label{3.16}
(\nabla_{X}Q)\xi = Q(\varphi +\varphi h)X - 2n\kappa(\varphi + \varphi h)X.
\end{eqnarray}
Taking the scalar product of (\ref{3.1}) with $\xi$ and using (\ref{3.16}) together with $Q\xi = 2n\kappa \xi$ (\ref{3.16}) gives
\begin{eqnarray}\label{3.17}
g(R(X,Y)D\lambda, \xi) &=& 2n\kappa[(X\lambda) - (Y\lambda)\eta(X)] + (\lambda + 1) g(Q\varphi X + \varphi QX, Y) \nonumber\\
&&+ (\lambda +1) g(Q\varphi hX + h\varphi QX,Y) - 4n\kappa (\lambda +1) g(\varphi X,Y)\nonumber\\
&&+ (Xf)\eta(Y) - (Yf)\eta(X).
\end{eqnarray}
Now, replacing  X by $\varphi X$, Y by $\varphi Y$ in (\ref{3.17}) and noting that $R(\varphi X,\varphi Y)\xi = 0$ (follows from (\ref{2.10})), we obtain
\begin{eqnarray*}
(\lambda + 1)[Q\varphi X + \varphi QX - \varphi QhX -hQ\varphi X - 4n\kappa \varphi X ]= 0.
\end{eqnarray*}
Since $\lambda$ is non-constant on $M$, the foregoing equation reduces to
\begin{eqnarray}\label{3.18}
Q\varphi X + \varphi QX - \varphi QhX -hQ\varphi X - 4n\kappa \varphi X =0
\end{eqnarray}
Substituting X by $\varphi X$ in (\ref{2.11}) gives
$$Q\varphi X = [2(n-1) - n\mu]\varphi X + [2(n-1) + \mu]h\varphi X.$$
On the other hand, the action of h on the forgoing equation and the use of (\ref{2.12}) provides
$$hQ\varphi X = [2(n-1) - n\mu]h\varphi X - (\kappa - 1)[2(n-1) + \mu]\varphi X.$$
Also operating (\ref{2.11}) by $\varphi$ yields
$$\varphi QX = [2(n-1) - n\mu]\varphi X + [2(n-1) + \mu]\varphi hX. $$
Taking $hX$ instead of $X$ and using (\ref{2.12}) the last equation reduces to
$$\varphi QhX = [2(n-1) - n\mu]\varphi hX - (\kappa - 1)[2(n-1) + \mu]\varphi X.$$
Next, we use the last four equations in (\ref{3.18}) to obtain
\begin{eqnarray}\label{3.19}
\kappa (\mu-2)=\mu(n+1).
\end{eqnarray}
Substituting  $\xi$ instead of X in (\ref{3.17}) and using $Q\xi = 2n\kappa\xi$ and (\ref{3.13}) we achieve
\begin{eqnarray}\label{3.20}
\kappa D\lambda + \mu hD\lambda - 2n\kappa((\xi\lambda)\xi - D\lambda) - (\kappa(\xi\lambda) + \xi f)\xi + D f=0.
\end{eqnarray}
Now, contracting (\ref{3.1}) over X and noting that the scalar curvature is constant we get
\begin{eqnarray}\label{3.21}
rD\lambda + 2nDf =0.
\end{eqnarray}
Moreover, using (\ref{3.21}) in (\ref{3.20}), we have
\begin{eqnarray}\label{3.22}
0 &=& 2n\kappa D\lambda + 2n\mu hD\lambda - 4n^2\kappa((\xi\lambda)\xi - D\lambda) \nonumber\\
&&+ 2n(\kappa(\xi\lambda) + \xi f)\xi + rD\lambda.
\end{eqnarray}
On the other hand, from equation (\ref{3.2}) and $Q\xi = 2n\kappa\xi$, we deduce
\begin{eqnarray}\label{3.23}
\nabla_{\xi}D\lambda = [2n\kappa(\lambda +1) - f]\xi.
\end{eqnarray}
Therefore, differentiating (\ref{3.22}) along $\xi$ and using (\ref{3.14}), (\ref{3.15}), (\ref{3.23}), we obtain
\begin{eqnarray}\label{3.24}
0 &=& (2n\kappa + 4n^2 \kappa - r)[2n\kappa(\lambda + 1) - f]\xi + 2n\mu^2 h\varphi D\lambda \nonumber\\
&&- 4n^2\kappa\xi(\xi\lambda)\xi - 2n\xi(\kappa(\xi\lambda) + \xi f)\xi.
\end{eqnarray}
Now, the action of $\varphi$ on the foregoing equation provides
$\mu^2hD\lambda =0.$
Further, operating the previous equation by $h$ and recalling (\ref{2.12}) it follows that
$$\mu^2(\kappa-1)\varphi^2D\lambda =0.$$
Since $M$ is non-Sasakian, we have
either (i) $\mu = 0$,  or (ii) $\varphi^2D\lambda = 0.$\\

\noindent
\textbf{Case(i):} In this case, it follows from (\ref{3.19}) that $\kappa = 0$. Hence $R(X,Y)\xi = 0$, and therefore $M^{2n+1}$ is locally flat in dimension $3$ and in higher dimensions it is locally isometric to the trivial bundle $E^{(n+1)}\times S^n(4)$ (see \cite{Blair}).\\

\noindent
\textbf{Case(ii):} This  case yields $ D\lambda = (\xi \lambda)\xi $. Differentiating this along an arbitrary vector field $ X $ together with (\ref{2.1}) entails that $ \nabla_{X}D\lambda = X(\xi \lambda)\xi - (\xi \lambda)(\varphi X - \varphi h X) $.
 Since $ g(\nabla_{X}D\lambda, Y) = g(\nabla_{Y}D\lambda, X) $, the foregoing equation shows $ X(\xi \lambda)\eta(Y) - Y(\xi \lambda )\eta(X) + (\xi \lambda) d\eta(X, Y) = 0 $.
Replacing X by $\varphi X$ and Y by $\varphi Y$ and noting that $ d\eta $ is non-zero for any contact metric structure, it follows that $ \xi \lambda = 0 $. Hence, $ D\lambda =0 $, i.e. $ \lambda $ is constant, which is a contradiction. This completes the proof.

\noindent\\

\textbf{Acknowledgments:}
The author D. S. Patra is financially supported by the Council of Scientific and Industrial Research, India (grant no. 17-06/2012(i)EU-V).

$^1$
Department of Mathematics, \\
Chandernagore College\\
Hooghly: 712 136 (W.B.), INDIA\\
E-mail: aghosh\_70@yahoo.com\\

$^2$
Department of Mathematics, \\
Jadavpur University,  \\
188. Raja S. C. Mullick Road,\\
Kolkata:700 032, INDIA \\
E-mail: dhritimath@gmail.com

\end{document}